\documentclass[12pt]{article}
\usepackage{amsmath}
\usepackage{epsf}
\usepackage{amssymb}
\usepackage{graphicx}
\usepackage{ifpdf}
\usepackage{epstopdf}
\usepackage{caption}
\usepackage{subfigure}
\newtheorem{theorem}{Theorem}
\newtheorem{lemma}{Lemma}
\newtheorem{corollary}{Corollary}
\newtheorem{observation}{Observation}
\newtheorem{example}{Example}
\newtheorem{remark}{Remark}
\newtheorem{acknowledgements}{Acknowledgement}

\begin{document}
\pagestyle{empty}
\renewcommand{\thefootnote}{\fnsymbol{footnote}}

\begin{center}
{\bf \Large Conflict-free connection number and independence number of a graph}\footnote{This
work was supported by Hunan Provincial Natural Science Foundation(No. 2018JJ2454) and Hunan Education Department Foundation(No. 18A382)}
\vskip 5mm

{ Jing Wang$^1$, Meng Ji$^2$\\[2mm]
 $1$ Department of Mathematics, Changsha University, Changsha 410003, China\\
 $2$ College of Mathematical Science, Tianjin Normal University, Tianjin, China}\\[6mm]
\end{center}
\date{}

\noindent{\bf Abstract} An edge-colored graph $G$ is conflict-free connected if any two of its vertices are connected by a path, which contains a color used on exactly one of its edges. The conflict-free connection number of a connected graph $G$, denoted by $cfc(G)$, is defined as the minimum number of colors that are required in order to make $G$ conflict-free connected. In this paper, we investigate the relation between the conflict-free connection number and the independence number of a graph. We firstly show that $cfc(G)\le \alpha(G)$ for any connected graph $G$, and an example is given showing that the bound is sharp. With this result, we prove that if $T$ is a tree with $\Delta(T)\ge \frac{\alpha(T)+2}{2}$, then $cfc(T)=\Delta(T)$.
{\bf MR(2000) Subject Classification} 05C15, 05C40

\section{Introduction}
\label{secintro}

All graphs considered here are simple, finite and undirected. An edge-coloring of a graph $G$ is {\it proper} if any two adjacent edges in this coloring receive different colors. If $G$ is colored with a proper coloring, then we say that $G$ is {\it properly colored}.

The rainbow connection number was introduced by Chartrand {\it et al.} \cite{Char2008}. An edge-colored graph $G$ is called {\it rainbow connected} if any two vertices are connected by a path whose edges have pairwise distinct colors. The {\it rainbow connection number} of a connected graph $G$, denoted by $rc(G)$, is the smallest number of colors that are needed in order to make $G$ rainbow connected. Chakraborty {\it et al.} \cite{Chak2011} showed that given a graph G, deciding if $rc(G)=2$ is NP-complete. Bounds for the rainbow connection number of a graph have also been studied in terms of other graph parameters, see \cite{dongli2016,Kam2016,lll2012,Kemn2013,lisun2017} and the references therein.

As an extension of proper colorings and motivated by rainbow connections of graphs, Andrews {\it et al.} \cite{Andrews2016} and independently Borozan {\it et al.} \cite{Borozen2012} introduced the concept of proper connection of graphs. An edge-colored graph $G$ is called {\it properly connected} if any two vertices are connected by a path which is properly colored. The {\it proper connection number} of a connected graph $G$, denoted by $pc(G)$, is the smallest number of colors that are needed in order to make $G$ properly connected. One can find many results on proper connection, see \cite{Brause2017,Huang2017,li2015,li2018} {\it et al}.

Very recently, inspired by rainbow connection colorings and proper connection colorings of graphs and by conflict-free colorings of graphs and hypergraphs \cite{chei2013,chei2011,Fabri2016,Pach2009}, Czap {\it et al.} \cite{czap} introduced the concept of conflict-free connection of graphs. An edge-colored graph $G$ is {\it conflict-free connected} if any two vertices are connected by a path, which contains at least one color used on exactly one of its edges. This path is called a {\it conflict-free path}, and this coloring is called a {\it conflict-free connection coloring} of $G$. The {\it conflict-free connection number} of a connected graph $G$, denoted by $cfc(G)$, is defined as the minimum number of colors that are required in order to make $G$ conflict-free connected.

An easy observation is that a rainbow edge-coloring of a connected graph $G$ is a trivial conflict-free connection coloring, while the other way around is not true in general. Moreover, all above mentioned three parameters of a graph $G$ with order $n$ are bounded by $n-1$, since one may color the edges of a given spanning tree of $G$ with distinct colors and color the remaining edges with already used colors.  There is an extensive research concerning on this topic, see \cite{chang2019,chaDAM,changji2019,deng,li2019,liwu2019}.

Recall that an {\it independent set} in a graph $G$ is a set of vertices no two of which are adjacent. The cardinality of a maximum independence set in $G$ is called the {\it independence number} of $G$ and is denoted by $\alpha(G)$. The observation follows immediately from the concept.

\begin{observation}\label{obser1}
Let $G$ be a connected graph of order $n$. Then $1\le \alpha(G)\le n-1$. Moreover, $\alpha(G)=1$ if and only if $G=K_n$, $\alpha(G)=n-1$ if and only if $G=K_{1,n-1}$.
\end{observation}

Dong and Li \cite{dongli2016} gave a relation between the rainbow connection number and the independence number of a graph, they showed that if $G$ is a connected graph without pendant vertices, then $rc(G)\le 2\alpha(G)-1$. Inspired by these results, we try to investigate the relation between the conflict-free connection number and the independence number of a graph and obtain our first main result.

\begin{theorem}\label{th1}
Let $G$ be a connected graph of order $n$. Then
\begin{equation*}
1\le cfc(G)\le \alpha(G)\le n-1.
\end{equation*}
Moreover, $cfc(G)=1$ if and only if $\alpha(G)=1$, $cfc(G)=n-1$ if and only if $\alpha(G)=n-1.$
\end{theorem}

Czap {\it et al.} \cite{czap} proved that 2-connected graphs have conflict-free connection number 2, while deciding the conflict-free connection number of graphs with cut-edges is very difficult, including trees. Chang {\it et al.} \cite{chaDAM} came up with a rapid approach to obtain the conflict-free connection number of a tree when its maximum degree is large. Motivated by these results, we find a method to determine the conflict-free connection number of a tree in terms of independence number and obtain our second main result.

\begin{theorem}\label{th2}
Let $T$ be a tree with $\Delta(T)\ge \frac{\alpha(T)+2}{2}$. Then $cfc(T)=\Delta(T)$.
\end{theorem}

We organize this paper as follows. Some useful preliminaries are presented in Section \ref{secpre}. Then, the proofs of Theorem \ref{th1} and Theorem \ref{th2} can be given in Section \ref{secth1} and Section \ref{secth2}, respectively.

We end this section with some terminology. Let $G=(V(G),E(G))$ be a graph with vertex set $V(G)$ and edge set $E(G)$. We use $d_G(v)$, $N_G(v)$ and $\Delta(G)$ to denote the {\it degree} of $v$ in $G$, the {\it set of neighbours} of $v$ in $G$ and the {\it maximum degree} of $G$, respectively. For $e\in E(G)$, we denote by $G\setminus e$ the graph obtained from $G$ by deleting $e$. An edge $e$ is said to be a {\it cut-edge} of $G$ if $c(G\setminus e)=c(G)+1$, where $c(G)$ is the number of components of $G$. Let $G$ and $F$ be two graphs, we use $F\subseteq G$ to denote that $F$ is a subgraph of $G$. For notation not explained here, readers are referred to \cite{Bondy2008}.

\section{Preliminaries}
\label{secpre}

This section is devoted to state several results which concerning on the conflict-free connection number of graphs. Czap {\it et al.} \cite{czap} showed that it is easy to obtain the conflict-free connection number for 2-connected graphs.

\begin{lemma}(\cite{czap})\label{le1}
If $G$ is a 2-connected and non-complete graph, then $cfc(G)=2$.
\end{lemma}

Chang {\it et al.} \cite{chaDAM} and independently Deng {\it et al.} \cite{deng} extended the result of Lemma \ref{le1} to 2-edge-connected graphs in the following.

\begin{lemma}(\cite{chaDAM}, \cite{deng})\label{le2}
Let $G$ be a non-complete 2-edge-connected graph, then $cfc(G)=2$.
\end{lemma}

Compared with 2-edge-connected graphs, the problem of determining the conflict-free connection number of graphs with cut-edges is very difficult. This fact arises many authors' attention to obtain lower or upper bounds of $cfc(G)$ for a connected graph. Chang {\it et al.} \cite{chaDAM} gave sharp lower and upper bound of $cfc(G)$ and characterized graphs $G$ for which $cfc(G)=1$ or $cfc(G)=n-1$.

\begin{lemma}(\cite{chaDAM})\label{le3}
Let $G$ be a connected graph of order $n$ $(n\ge 2)$. Then $1\le cfc(G)\le n-1.$ Moreover, $cfc(G)=1$ if and only if $G=K_n$, $cfc(G)=n-1$ if and only if $G=K_{1,n-1}$.
\end{lemma}

A {\it block} of a graph $G$ is a maximal connected subgraph of $G$ that has no cut-vertex. If $G$ itself is connected and has no cut-vertex, then $G$ is a block. An edge is a block if and only if it is a cut-edge. A block consisting of a cut-edge is called trivial. Note that any nontrivial block is 2-connected.

Let $C(G)$ be the subgraph of $G$ induced on the set of cut-edges of $G$, and let $h(G)=\max\{cfc(T): T {\rm ~ is~ a ~component ~of}~ C(G)\}$.

\begin{lemma}(\cite{czap})\label{le4}
If $G$ is a connected graph with cut-edges, then $h(G)\le cfc(G) \le h(G)+1$. Moreover, these bounds are tight.
\end{lemma}

Chang {\it et al.} \cite{chaDAM} gave a sufficient condition such that the lower bound in Lemma \ref{le4} is sharp for $h(G)\ge 2.$

\begin{lemma}(\cite{chaDAM})\label{le5}
Let $G$ be a connected graph with $h(G)\ge 2$. If there exists a unique component $T$ of $C(G)$ satisfying $(i)$ $cfc(T)=h(G)$, $(ii)$ $T$ has an optimal conflict-free connection coloring which contains a color used on exactly one edge of $T$, then $cfc(G)=h(G)$.
\end{lemma}

It is seen from Lemma \ref{le4} that, to determine the conflict-free connection number of graphs relies on the conflict-free connection number of trees, with an error of only one. Thus, determining the conflict-free connection number of trees is of great importance. Here we list some known results concentrating on the conflict-free connection number of trees.

\begin{lemma}(\cite{czap})\label{le6}
If $P_n$ is a path on $n$ edges, then $cfc(P_n)=\big\lceil \log_2(n+1)\big\rceil$.
\end{lemma}

\begin{lemma}(\cite{czap})\label{le7}
If $T$ is an $n-$vertex tree of maximum degree $\Delta(T)\ge 3$ and diameter $d(T)$, then
\begin{equation*}
\max\big\{\Delta(T), \log_2d(T)\big\}\le cfc(T)\le \frac{\big(\Delta(T)-2\big)\cdot \log_2n}{\log_2\Delta(T)-1}.
\end{equation*}
\end{lemma}

The following result in \cite{chaDAM} indicates that when the maximum degree of a tree is large, the conflict-free connection number is immediately determined by its maximum degree.

\begin{lemma}(\cite{chaDAM})\label{le8}
Let $T$ be a tree of order $n$, and $t$ be a positive integer such that $n\ge 2t+2$. Then $cfc(T)=n-t$ if and only if $\Delta(T)=n-t.$
\end{lemma}

We end this section with the following lemma, which is no more than an observation.

\begin{lemma}\label{le9}
Let $T_1$ and $T_2$ be two trees such that $T_1\subseteq T_2$. Then $cfc(T_1)\le cfc(T_2)$.
\end{lemma}

\section{The proof of Theorem \ref{th1}}
\label{secth1}

\noindent{\bf Proof of Theorem \ref{th1}.} \, By Observation \ref{obser1} and Lemma \ref{le3}, it suffices to prove that $cfc(G)\le \alpha(G)$ for a non-complete graph $G$. Our main strategy is by induction on the number of cut-edges in $G$. For simplicity, set $k:=|E(C(G))|$.

Since $cfc(G)=2$ and $\alpha(G)\ge 2$ for a non-complete 2-edge-connected graph $G$, we get that $cfc(G)\le \alpha(G)$ when $k=0$ by Lemma \ref{le2} and Observation \ref{obser1}. Assume that the statement holds for any graph with $\le k-1$ cut-edges, and let $G$ be a graph with $k$ cut-edges. We distinguish two cases.

\vskip 1mm
{\bf Case 1.} There exists a cut-edge, say $e$, such that each component of $G\setminus e$ is a subgraph of order greater than 1.

W.l.o.g., let $e=u_1u_2$ and let $G_1$ and $G_2$ be two components of $G\setminus e$ with $u_i\in V(G_i)$, $i\in\{1,2\}$.

For $i\in\{1,2\}$, it is seen that $|V(G_i)|\ge 2$, and that the number of cut-edges in $G_i$ must be no more than $k-1$. By induction hypothesis, we have
\begin{equation*}
cfc(G_i)\le \alpha(G_i) ~~~ {\rm for}~~ i\in\{1,2\}.
\end{equation*}

W.l.o.g., assume that $cfc(G_2)\le cfc(G_1)$. Let $S_1$ be a maximum independent set in $G_1$. Moreover, since $|V(G_2)|\ge 2$, there must exist a vertex, say $z$, such that $z\in V(G_2)\setminus \{u_2\}$. Note that $z$ is not adjacent to vertices in $G_1$, then $S_1\cup\{z\}$ is an independent set in $G$ whose cardinality is
\begin{equation*}
\alpha(G_1)+1\le \alpha(G).
\end{equation*}

Now, we are able to assign $cfc(G_1)+1$ colors to all the edges of $G$ in order to make $G$ is conflict-free connected: first we color each component of $G\setminus e$ with at most $cfc(G_1)$ colors, next we color the edge $e$ with a fresh color. We only need to prove that any pair of distinct vertices $x$ and $y$ of $G$ are connected by a conflict-free path. If the vertices $x$ and $y$ are from the same component of $G\setminus e$, then such a path exists. If they are in different components of $G\setminus e$, then there is a $x-y$ path through the edge $e$ with a unique color.

The analyses above imply that
\begin{equation*}
cfc(G)\le cfc(G_1)+1\le \alpha(G_1)+1\le \alpha(G).
\end{equation*}

\vskip 1mm
{\bf Case 2.} Each cut-edge is a pendant edge.

Thus, each component of $C(G)$ is a complete bipartite graph $K_{1,r}$ where $1\le r\le n-1.$ Let $\widetilde{G}$ be the graph obtained from $G$ by deleting all the pendant vertices. Note that $|V(\widetilde{G})|\ne 2$, otherwise $\widetilde{G}$ is a non-pendant cut-edge in $G$, a contradiction.

\vskip 1mm
{\bf Subcase 2.1.} $|V(\widetilde{G})|=1$.

That means $G=K_{1,n-1}$. By Observation \ref{obser1} and Lemma \ref{le3}, $cfc(G)=n-1=\alpha(G).$

\vskip 1mm
{\bf Subcase 2.2.} $|V(\widetilde{G})|\ge 3$.

W.l.o.g., let $v$ be a vertex of $C(G)$ such that
$$d_{C(G)}(v)=\max\big\{d_{C(G)}(x): x\in V(C(G))\big\}.$$
For simplicity, setting $t:=d_{C(G)}(v)$ and let $y_1,\cdots, y_t$ be pendant vertices adjacent to $v$ in $G$. Thus,
\begin{equation}\label{eq00}
h(G)=cfc(K_{1,t})=t.
\end{equation}

Since $|V(\widetilde{G})|\ge 3$, we can choose a vertex, say $z$, such that $z\in V(\widetilde{G})\setminus \{v\}$. Note that $\{z, y_1,\cdots, y_t\}$ is an independent set in $G$ with cardinality $t+1$, obviously
\begin{equation}\label{eq1}
t+1\le \alpha(G).
\end{equation}

Therefore, Lemma \ref{le4} together with Eq.(\ref{eq00}) and Eq.(\ref{eq1}) yield
\begin{equation*}
cfc(G)\le h(G)+1=t+1\le \alpha(G).
\end{equation*}

Thus, $1\le cfc(G)\le \alpha(G) \le n-1$ for a connected graph $G$.

Moreover, Observation \ref{obser1} together with Lemma \ref{le3} imply that $cfc(G)=1$ if and only if $\alpha(G)=1$, and that $cfc(G)=n-1$ if and only if $\alpha(G)=n-1.$ \hfill$\Box$

\vskip 1mm
By Theorem \ref{th1}, it is easy to obtain the conflict-free connection number of a graph whose independence number is 2.

\begin{corollary}\label{co1}
Let $G$ be a connected graph with $\alpha(G)=2$. Then $cfc(G)=2.$
\end{corollary}

By Theorem \ref{th1} and Lemma \ref{le7}, we can give an upper bound on the conflict-free connection number of trees. Moreover, a sufficient condition for which the conflict-free connection number of a tree equals to its maximum degree is obtained.

\begin{corollary}\label{co2}
Let $T$ be a tree. Then $\Delta(T)\le cfc(T)\le \alpha(T).$ Moreover, if $\Delta(T)=\alpha(T)$, then $cfc(T)=\Delta(T).$
\end{corollary}

At the end of this section, an example is given showing that there exists non-complete graph whose conflict-free connection number can be any integer no more than its independence number. Thus the bound $cfc(G)\le \alpha(G)$ in Theorem \ref{th1} is tight.

\begin{example}\label{ex1}
Let $l,k$ be integers such that $3\le l\le n-2$ and that $2\le k\le l$. There exists a graph $G_{l,k}$ of order $n$ for which $\alpha(G_{l,k})=l$ and $cfc(G_{l,k})=k.$
\end{example}

{\bf Proof.} We will construct the desired graph $G_{l,k}$ by considering two cases: $k=l$ or $k<l.$

When $k=l,$ let $G_{l,l}$ be a graph obtained by identifying a leaf vertex of $K_{1,l}$ with a vertex of the complete graph $K_{n-l}$. It is seen that $\alpha(G_{l,l})=l=cfc(G_{l,l})$.

When $k<l,$ we construct $G_{l,k}$ with vertex set
\begin{equation*}
V(G_{l,k})=\{w,v, u_1,\cdots, u_{n-2}\}
\end{equation*}
and edge set
\begin{align*}
E(G_{l,k})=&  \{wu_i: 1\le i\le n-2\}\,\cup \,\{vu_i: k+1\le i\le n-2\}\\
& \cup \,\{wv\}\,\cup \,\{u_iu_j: l\le i\ne j\le n-2\}.
\end{align*}

Note that the subgraph induced on vertices $\{v, w, u_l, \cdots, u_{n-2}\}$ is a clique on $n-l+1$ vertices. We can get that $\alpha(G_{l,k})=l$ since $\{u_1, \cdots, u_k, \cdots, u_l\}$ is a maximum independent set in $G_{l,k}.$ Moreover, the subgraph induced on $\{w,u_1, \cdots, u_k\}$ is the unique component of $C(G_{l,k})$, thus $cfc(G_{l,k})=k$ by Lemma \ref{le3} and Lemma \ref{le5}. \hfill$\Box$

\section{The proof of Theorem \ref{th2}}
\label{secth2}

We firstly give some results on the conflict-free connection number of certain trees, which will be useful in the later discussions.

\begin{figure}[htbp]
\begin{minipage}[t]{0.5\linewidth}
\resizebox{0.8\textwidth}{!} {\includegraphics{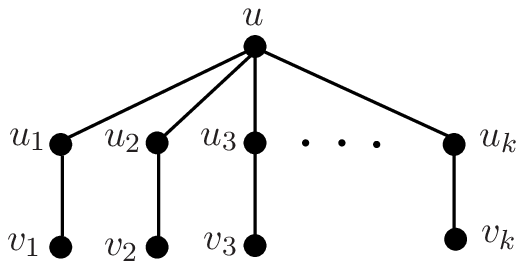}}
\caption{\small The graph $H_{k}.$ } \label{fig1}
\end{minipage}
\begin{minipage}[t]{0.5\linewidth}
\resizebox{0.8\textwidth}{!} {\includegraphics{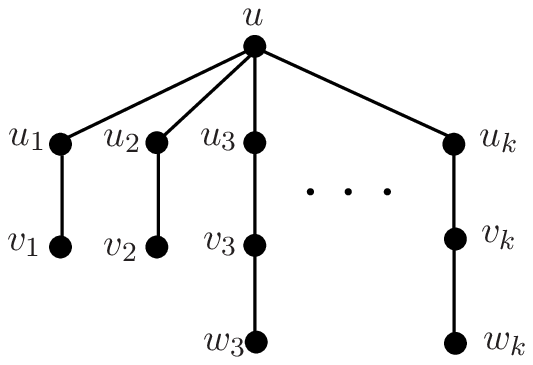}}
\caption{\small The graph $Q_k.$} \label{fig2}
\end{minipage}
\end{figure}

\begin{lemma}\label{le010}
Define $H_k$ $(k\ge 3)$ be a tree obtained by subdividing each edge of the complete bipartite graph $K_{1,k}$ to a path of length two, see Figure \ref{fig1}. Then $cfc(H_k)=k.$
\end{lemma}

{\bf Proof.} By the definition of $H_k$ and by Lemma \ref{le7}, we have $cfc(H_k)\ge \Delta(H_k)=k.$ To complete the proof, we only need to assign a conflict-free connected coloring $c:\, E(H_k)\rightarrow [k]$ as follows
\[
c(e)= \left\{
\begin{array}{ll}
 i, \quad & {\rm if}~ e=uu_i, \;\; 1\le i\le k; \\
 k, \quad &  {\rm if}~ e=u_1v_1;\\
 i-1, \quad & {\rm if}~ e=u_iv_i, \;\; 2\le i\le k.
\end{array}
\right.
\]
It is not difficult to check that $c$ is a conflict-free connected coloring of $H_k$, thus $cfc(H_k)\le k.$ The proof is done. \hfill$\Box$

\begin{lemma}\label{le011}
Define $Q_k$ $(k\ge 3)$ be a tree obtained from $H_{k}$ by adding a pendant edge to each of the $k-2$ leaf vertex of $H_{k}$, see Figure \ref{fig2}. Then $cfc(Q_k)=k.$
\end{lemma}

{\bf Proof.} By Lemma \ref{le7}, we only need to assign a conflict-free connected coloring $c:\, E(Q_k)\rightarrow [k]$ as follows
\[
c(e)= \left\{
\begin{array}{ll}
 i, \quad &  {\rm if}~ e=uu_i, \;\; 1\le i\le k; \\
 k, \quad & {\rm if}~  e=u_1v_1;\\
 i-1, \quad & {\rm if}~ e=u_iv_i, \;\; 2\le i\le k;\\
 1, \quad & {\rm if}~ e=v_iw_i, \;\; 3\le i\le k.
\end{array}
\right.
\]\hfill$\Box$

Now we are in a position to prove Theorem \ref{th2}.

\vskip 1mm
{\bf Proof of Theorem \ref{th2}.} We will prove the theorem by induction on $k:=\Delta(T)$.

Since a tree $T$ satisfying $\Delta(T)=2$ and $\Delta(T)\ge (\alpha(T)+2)/2$ is the path $P_2$ or $P_3$, by Lemma \ref{le6}, the theorem holds when $k=2$. Assume that the result is true for any tree $T'$ with $\Delta(T')\le k-1$ and $\Delta(T')\ge (\alpha(T')+2)/2$. Now consider a tree $T$ with $\Delta(T)=k$ $(k\ge 3)$ and $\Delta(T)\ge (\alpha(T)+2)/2$.

Let $u$ be a vertex of $T$ such that $d_T(u)=\Delta(T)$ and let $N_T(u)=\{u_1, \cdots, u_k\}$. Firstly, we claim that

\vskip 1mm
{\bf Claim 1.} If there exists a vertex $w\ne u$ such that $d_T(w)=\Delta(T)$, then $w\in N_T(u)$.

\vskip 1mm
{\bf Proof of Claim 1.} Suppose to contrary that there is a vertex $w\notin N_T(u)$ such that $d_T(w)=\Delta(T)$. Then $N_T(u)\cup N_T(w)$ is an independent set in $T$, whose cardinality is at least
\begin{equation*}
2k-1=2\Delta(T)-1\ge \alpha(T)+1,
\end{equation*}
the last inequality holds since $\Delta(T)\ge (\alpha(T)+2)/2$. A contradiction. \hfill$\Box$

\vskip 1mm
It is inferred from the proof of Claim 1 that, in $T$, the vertices of maximum degree must be adjacent to each other. Since $T$ is a tree, there is at most one vertex of $\{u_1, \cdots, u_k\}$ can be of maximum degree. W.l.o.g., let
\begin{equation*}
d_T(u_1)=\max \{d_T(u_i): 1\le i\le k\}.
\end{equation*}

Therefore, we claim that

\vskip 1mm
{\bf Claim 2.} For any vertex $x\in V(T)\setminus\{u,u_1\}$, it has $d_T(x)\le \Delta(T)-1$.

\vskip 1mm

For $1\le i\le k,$ let $T_{i1}$ and $T_{i2}$ be two components of $T\setminus uu_i$ where $u_i\in V(T_{i2})$. We discuss three cases.

\vskip 1mm
{\bf Case 1.} $d_{T}(u_1)=1.$
\vskip 1mm

Then $T$ is the graph $K_{1,k}$. By Lemma \ref{le3}, $cfc(T)=k=\Delta(T)$.

\vskip 1mm
{\bf Case 2.} $d_{T}(u_1)=2$.

\vskip 1mm
{\bf Subcase 2.1.} $|E(T_{i2})|\le 1$ for each $1\le i\le k.$
\vskip 1mm

Then $T$ is a subgraph of $H_k$ which is defined in Lemma \ref{le010}. By Lemma \ref{le9} and Lemma \ref{le010}, we have $cfc(T)\le cfc(H_k)=k$. On the other hand, by Lemma \ref{le7}, $cfc(T)\ge \Delta(T)=k.$ Thus, $cfc(T)=k=\Delta(T).$

\vskip 1mm
{\bf Subcase 2.2.} $|E(T_{i2})|\le 2$ for each $1\le i\le k,$ and there exists an integer $i$ such that $|E(T_{i2})|=2$.
\vskip 1mm

Note that $T_{i2}$ is a path $P_2$ when $|E(T_{i2})|=2$. W.l.o.g., let $w_i$ be an end vertex other than $u_i$ in $T_{i2}$. Since $d_{T}(u_1)=2$, there are at most $k-2$ integers $i$ such that $T_{i2}$ is a path $P_2$, otherwise $S:=\{u_1, \cdots, u_k\}\cup \big(\bigcup_i\{w_i\}\big)$ is an independent set in $T$, moreover,
\begin{equation*}
|S|\ge k+k-1=2\Delta(T)-1\ge \alpha(T)+1,
\end{equation*}
a contradiction.

Therefore, $T$ is a subgraph of $Q_k$ which is defined in Lemma \ref{le011}. By Lemmas \ref{le7}, \ref{le9} and \ref{le011}, we have $cfc(T)=k=\Delta(T).$

\vskip 1mm
{\bf Subcase 2.3.} There exists an integer $i$ $(1\le i\le k)$ such that $|E(T_{i2})|\ge 3$.

\vskip 1mm
Let $e=uu_i$. Recall that $T_{i1}$ and $T_{i2}$ are two components of $T\setminus e$ with $u_i\in V(T_{i2})$. Since $d_{T}(u_1)=2<\Delta(T)$, thus $\Delta(T_{i1})=\Delta(T)-1$ and $\Delta(T_{i2})\le\Delta(T)-1$ by Claim 2.

Firstly, we try to obtain the conflict-free connection number of $T_{i1}$. Let $S_1$ be a maximum independent set in $T_{i1}$. Since $|E(T_{i2})|\ge 3$, we always can choose at least two non-adjacent vertices, say $x$ and $y$, from $V(T_{i2})\setminus \{u_i\}$, such that $S:=S_1\cup \{x,y\}$ is an independent set in $T$. That means $|S|=\alpha(T_{i1})+2\le \alpha(T)$, therefore,
\begin{equation*}
\Delta(T_{i1})=\Delta(T)-1\ge \frac{\alpha(T)}{2}\ge \frac{\alpha(T_{i1})+2}{2},
\end{equation*}
the above first inequality holds since $\Delta(T)\ge \frac{\alpha(T)+2}{2}$. By induction hypothesis, we have
\begin{equation}\label{eq2}
cfc(T_{i1})=\Delta(T_{i1})=\Delta(T)-1.
\end{equation}

Next, we consider the conflict-free connection number of $T_{i2}$. Firstly, we claim that

\vskip 1mm
{\bf Claim 3.} $\alpha(T_{i2})\le \alpha(T)-\Delta(T)+1.$

\vskip 1mm
{\bf Proof of Claim 3.} Suppose to contrary that $\alpha(T_{i2})> \alpha(T)-\Delta(T)+1$. Let $S_2$ be a maximum independent set in $T_{i2}$, obviously, $S':=S_2\cup \big(\bigcup_{j\ne i}\{u_j\}\big)$ is an independent set in $T$ with cardinality
\begin{equation*}
|S'|=|S_2|+k-1>\alpha(T)-\Delta(T)+1+k-1=\alpha(T),
\end{equation*}
a contradiction.  \hfill$\Box$

By Theorem \ref{th1} and Claim 3, we have
\begin{equation}\label{eq3}
cfc(T_{i2})\le \alpha(T_{i2})\le \alpha(T)-\Delta(T)+1\le \Delta(T)-1,
\end{equation}
the last inequality holds since $\Delta(T)\ge \frac{\alpha(T)+2}{2}$.

By Eq.(\ref{eq2}) and Eq.(\ref{eq3}), we are now able to assign $\Delta(T)$ colors to all the edges of $T$ in order to make $T$ is conflict-free connected: firstly we color $T_{i1}$ and $T_{i2}$ with at most $\Delta(T)-1$ colors, next we color the edge $e=uu_i$ with a fresh color. Therefore, $cfc(T)\le \Delta(T).$ Combined this conclusion with Lemma \ref{le7}, we have $cfc(T)=\Delta(T).$

\vskip 1mm
{\bf Case 3.} $d_{T}(u_1)\ge 3$.

\vskip 1mm

Let $e=uu_1$. Recall that $T_{11}$ and $T_{12}$ are two components of $T\setminus e$ with $u_1\in V(T_{12})$. By Claim 2, $\Delta(T_{11})=\Delta(T)-1$ and $\Delta(T_{12})\le\Delta(T)-1$.

Using similar discussions in Subcase 2.3, we can get that $cfc(T_{11})=\Delta(T)-1$ and that $cfc(T_{12})\le\Delta(T)-1$, moreover, $cfc(T)=\Delta(T).$

The proof is completed. \hfill$\Box$

\vskip 1mm
\begin{remark}\label{re1}
The sharpness example for Theorem \ref{th2} is given as follows. Let $T$ be a tree obtained from two copies of $K_{1,k-1}$ with $k\ge 3$ by identifying a leaf vertex in one copy with a leaf vertex in the other copy. It is seen that $\alpha(T)=2k-3$ and that $\Delta(T)=k-1=\frac{\alpha(T)+1}{2}$. Furthermore, Theorem 5.5 in \cite{changji2019} showed that $cfc(T)=k$. Thus  $cfc(T)>\Delta(T)$.
\end{remark}

\begin{figure}[htbp]
\centering
\resizebox{0.4\textwidth}{!} {\includegraphics{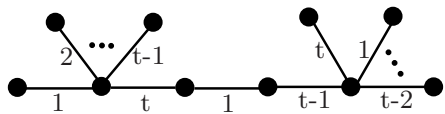}}
\caption{\small The graph $G.$ } \label{fig3}
\end{figure}

\begin{remark}\label{re2}
Theorem \ref{th2} gives a sufficient condition for the conflict-free connection number of a tree equals to its maximum degree. However, this condition is not necessary. Define $G$ to be a tree obtained from two copies of $K_{1,k}$ $(k\ge 3)$ by adding an edge joining a leaf vertex in one copy to a leaf vertex in the other copy. Figure \ref{fig3} illustrates that we can assign $k=\Delta(G)$ colors to all the edges of $G$ in order to make it conflict-free connected, thus $cfc(G)=\Delta(G)$. On the other hand, we can testify that $\alpha(G)=2k-1$ and thus $\Delta(G)<\frac{\alpha(G)+2}{2}.$
\end{remark}

\begin{acknowledgements}
This work was done during the first author was visiting Nankai University. The authors would like to express their sincere thanks to Professor Xueliang Li for his helpful suggestions.
\end{acknowledgements}


\end{document}